\input amstex
\documentstyle{amsppt}
\document
\topmatter
\title
K\"ahler manifolds with quasi-constant  holomorphic curvature.
\endtitle
\author
W\l odzimierz Jelonek
\endauthor

\abstract{The aim of this paper is to classify compact  K\"ahler
manifolds with quasi-constant holomorphic sectional curvature. }
\thanks{MS Classification: 53C55,53C25. Key words and phrases:
K\"ahler manifold, holomorphic sectional curvature, quasi constant
holomorphic sectional curvature, special K\"ahler-Ricci
potential}\endthanks
 \endabstract
\endtopmatter
\define\G{\Gamma}
\define\DE{\Cal D^{\perp}}
\define\e{\epsilon}
\define\n{\nabla}
\define\om{\omega}
\define\w{\wedge}
\define\k{\diamondsuit}
\define\th{\theta}
\define\p{\partial}
\define\a{\alpha}
\define\be{\beta}
\define\g{\gamma}
\define\lb{\lambda}

\define\1{D_{\lb}}
\define\2{D_{\mu}}
\define\0{\Omega}

\define\De{\Cal D}

\define\dl{\delta}

\define\m{(M,g,J)}
\define \E{\Cal E}
\bigskip
{\bf 0. Introduction. } The aim of the present paper is to
classify compact, simply connected K\"ahler manifolds $(M,g,J)$
admitting a global, $2$-dimensional, $J$-invariant  distribution
$\De$ having the following property:  The holomorphic curvature
$K(\pi)=R(X,JX,JX,X)$ of any $J$-invariant $2$-plane $\pi\subset
T_xM$, where $X\in \pi$ and $g(X,X)=1$, depends only on the point
$x$ and the number $|X_{\De}|=\sqrt{g(X_{\De},X_{\De})}$, where
$X_{\De}$ is an orthogonal projection of $X$ on $\De$. In this
case  we have
$$R(X,JX,JX,X)=\phi(x,|X_{\De}|)$$ where $\phi(x,t)=a(x)+b(x)t^2+c(x)t^4$ and
 $a,b,c$ are smooth functions on $M$. Also $R=a\Pi+b\Phi+c\Psi$
 for certain curvature tensors $\Pi,\Phi,\Psi\in \bigotimes^4\frak X^*(M)$
  of K\"ahler type. The investigation of such manifolds, called QCH K\"ahler manifolds,  was
started by G. Ganchev and V. Mihova in [G-M-1],[G-M-2]. In our
paper we shall use their local results to obtain a global
classification of such manifolds under the assumption that $\dim
M=2n\ge 6$.
 By $\E$ we shall denote the $(\dim M-2)$-dimensional distribution
 which is the orthogonal complement of $\De$ in $TM$. If $X$ is a local unit section of $\De$
 then  $\{X,JX\}$
 is a local orthonormal basis of $\De$ and the function $\kappa=
 \sqrt{(div_{\E}X)^2+(div_{\E}JX)^2}$ does not depend on the
 choice of $X$. We classify  compact, simply connected QCH K\"ahler manifolds
 satisfying the conditions
 int
$B=\emptyset$ and $U\ne \emptyset$ where $B=\{x\in U:
b(x)=0\},U=\{x\in M: \kappa(x)\ne 0\}$. First we  show that $\m$
admits a global holomorphic Killing vector field with a Killing
potential, which is a special K\"ahler-Ricci potential. Next we
use the results of Derdzi\'nski and Maschler [D-M-1], who
classified compact K\"ahler manifolds admitting special
K\"ahler-Ricci potentials. As a corollary we prove that every
compact, simply connected QCH K\"ahler manifold with $\kappa\ne 0$
and an analytic Riemannian metric $g$ is a holomorphic
$\Bbb{CP}^1$-bundle over $\Bbb{CP}^{n-1}$ (with $\De$ being an
integrable distribution whose leaves are the  $\Bbb{CP}^1$ fibers
of the bundle) or it is $\Bbb{CP}^n$ with a metric of constant
holomorphic sectional curvature ( in this case $\De$ is any
$J$-invariant 2-dimensional distribution on $\Bbb{CP}^n$ with
$\kappa\ne 0$; however, such distributions may not exist).
\bigskip
{\bf 1. Special frame.} Let $\m$ be a $2n$-dimensional K\"ahler
manifold with a $2$-dimensional $J$-invariant distribution $\De$.
Let $\frak X(M)$ denote the algebra of all differentiable vector
fields on $M$ and  $\G(\De)$ denote the set of local sections of
the distribution $\De$. If $X\in\frak X(M)$ then by $X^{\flat}$ we
shall denote the 1-form $\phi\in \frak X^*(M)$ dual to $X$ with
respect to $g$, i.e. $\phi(Y)=X^{\flat}(Y)=g(X,Y)$. By $\0$ we
shall denote the K\"ahler form of $\m$ i.e. $\0(X,Y)=g(JX,Y)$. Let
us denote by $\E$ the distribution $\DE$, which is a
$2(n-1)$-dimensional, $J$-invariant distribution. By $h,m$
respectively we shall denote the tensors $h=g\circ (p_{\De}\times
p_{\De}),m=g\circ (p_{\E}\times p_{\E})$, where $p_{\De},p_{\E}$
are the orthogonal projections on $\De,\E$ respectively. It
follows that $g=h+m$. By $\om$ we shall denote the K\"ahler form
of $\De$ i.e. $\om(X,Y)=h(JX,Y)$ and by $\0_m$ the K\"ahler form
of $\E$ i.e. $\0_m(X,Y)=m(JX,Y)$.    For any local section
$X\in\G(\De)$ we define $div_{\E}X=tr_m\n
X^{\flat}=m^{ij}\n_{e_i}X^{\flat}(e_j)$ where
$\{e_1,e_2,...,e_{2(n-1)}\}$ is any basis of $\E$ and $[m^{ij}]$
is a matrix inverse to $[m_{ij}]$, where $m_{ij}=m(e_i,e_j)$. Note
that if $f\in C^{\infty}(M)$ then $div_{\E}(f X)=fdiv_{\E}X$ in
the case  $X\in \G(\De)$. Let $\xi\in\G(\De)$ be a unit local
section of $\De$. Then $\{\xi,J\xi\}$ is an orthonormal basis of
$\De$. Let $\th(X)=g(\xi,X)$ and $J\th=-\th\circ J$ which means
that $J\th(X)=g(J\xi,X)$.  Let us denote by $\kappa$ the function
$$\kappa=\sqrt{(div_{\E}\xi)^2+(div_{\E}J\xi)^2}.\tag 1.1$$
The function $\kappa$ does not depend on the choice of a section
$\xi$. In fact, if $\xi'=a\xi+b J\xi$, where $a,b\in
C^{\infty}(dom \xi)$ and $a^2+b^2=1$ is another unit section of
$\De$, then $J\xi'=-b\xi+aJ\xi$ and
$$\gather
(div_{\E}\xi')^2+(div_{\E}J\xi')^2=(a div_{\E}\xi+b
div_{\E}J\xi)^2+(-b div_{\E}\xi+a div_{\E}J\xi)^2\tag 1.2
\\=(div_{\E}\xi)^2+(div_{\E}J\xi)^2.\endgather$$ Hence $\kappa$ is
a well defined, continuous function on $M$, which is smooth in the
open set $U=\{x:\kappa(x)\ne 0\}$. We shall now show that on $U$
there is a smooth, global unit  section $\xi\in \G(U, \De)$
defined uniquely up to  a sign such that $div_{\E}J\xi=0$. Namely,
if $\xi'$ is a local unit section of $\G(U,\De)$ then then the
section $\xi=\frac 1{\kappa}
((div_{\E}\xi')\xi'+(div_{\E}J\xi)J\xi')$ satisfies
$div_{\E}J\xi=\frac1{\kappa}((div_{\E}\xi')(div_{\E}J\xi')-(div_{\E}J\xi')(div_{\E}\xi'))=0$
and does not depend on the choice of $\xi'$. On the other hand it
is clear that the only other such smooth section is $-\xi$. The
section $\xi$ constructed above and defined on $U\subset M$ we
shall call the principal section of $\De$. Note that
$div_{\E}\xi=\kappa$.
\medskip
{\bf 2. Curvature tensor of a QCH K\"ahler manifold.} We shall
recall some results from [G-M-1]. Let
$R(X,Y)Z=([\n_X,\n_Y]-\n_{[X,Y]})Z$ and let us write $$R(X
,Y,Z,W)=g(R(X,Y)Z,W).$$ If $R$ is the curvature tensor of a QCH
K\"ahler manifold $\m$, then there exist functions $a,b,c\in
C^{\infty}(M)$ such that
$$R=a\Pi+b\Phi+c\Psi,\tag 2.1$$
where $\Pi$ is the standard K\"ahler tensor of constant
holomorphic curvature i.e.
$$\gather \Pi(X,Y,Z,U)=\frac14(g(Y,Z)g(X,U)-g(X,Z)g(Y,U)\tag 2.2\\+g(JY,Z)g(JX,U)-g(JX,Z)g(JY,U)-2g(JX,Y)g(JZ,U)),\endgather $$
the tensor $\Phi$ is defined by the following relation
$$\gather \Phi(X,Y,Z,U)=\frac18(g(Y,Z)h(X,U)-g(X,Z)h(Y,U)\tag 2.3\\+g(X,U)h(Y,Z)-g(Y,U)h(X,Z)
+g(JY,Z)\om(X,U)\\-g(JX,Z)\om(Y,U)+g(JX,U)\om(Y,Z)-g(JY,U)\om(X,Z)\\
-2g(JX,Y)\om(Z,U)-2g(JZ,U)\om(X,Y)),\endgather$$ and finally
$$\Psi(X,Y,Z,U)=-\om(X,Y)\om(Z,U)=-(\om\otimes\om)(X,Y,Z,U).\tag 2.4$$

 Let $V=
 (V,g,J)$ be a real $2n$ dimensional vector space with
complex structure $J$ which is skew-symmetric with respect to the
scalar product $g$ on $V$.  Let assume further that  $V= D\oplus
E$ where $D$ is a 2-dimensional, $J$-invariant subspace of $V$,
$E$ denotes its orthogonal complement in $V$. Note that  the
tensors $\Pi,\Phi,\Psi$ given above are of K\"ahler type.  It is
easy to check that for a unit vector $X\in V$
$\Pi(X,JX,JX,X)=1,\Phi(X,JX,JX,X)=|X_{D}|^2,\Psi(X,JX,JX,X)=|X_{D}|^4$,
where $X_D$ means an orthogonal projection of a vector $X$ on the
subspace $D$ and $|X|=\sqrt{g(X,X)}$. It follows that for a tensor
$(2.1)$ defined on $V$ we have
$$R(X,JX,JX,X)=\phi(|X_D|)$$ where $\phi(t)=a+bt^2+ct^4$.

Now let us assume that $\m$ is a QCH K\"ahler manifold of
dimension $2n\ge 6$ and let $X$ be a local  unit section of $\De$
and $\eta(Z)=g(X,Z)$. Let us define two 1-forms $\e,\e^*$ by the
formulas: $\e(Z)=g(p_{\E}(\n_XX),Z)=g(\n_XX,Z)-pJ\eta(Z)$,
$\e^*(Z)=g(p_{\E}(\n_{JX}JX),Z)=g(\n_{JX}JX,Z)-p^*\eta(Z)$ where
$p=g(\n_XX,JX)$,$p^*=g(\n_{JX}JX,X)$ and $p_{\E}$ denotes the
orthogonal projection on $\E$. Note that the distribution $\De$ is
integrable if and only if $\e+\e^*=0$ (see [G-M-1], Lemma 3.3). In
fact for $Z\in\G(\E)$ we have
$$\gather g([X,JX],Z)=g(\n_XJX-\n_{JX}X,Z)=g(J\n_XJX-J\n_{JX}X,JZ)\\
=-g(\n_XX+\n_{JX}JX,JZ)=-(\e(JZ)+\e^*(JZ)).\endgather$$
 Let $\{Z_{\lb}\}$ be any complex basis of the complex
subbundle $\E^{1,0}$ of the complex tangent bundle $T^cM=\Bbb
C\otimes TM$. We also write $Z_{\bar\lb}=\bar{Z_{\lb}}$. Then the
Bianchi identity for the tensor $R$ of the form $(2.1)$ gives the
following relations (see Theorem 3.5 in [G-M-1]) :
$$\gather
\n a=\frac{b div_{\E}X}{2(n-1)}X+\frac{b
div_{\E}JX}{2(n-1)}JX,\tag 2.5\\
\n b=\frac{(b+4c)div_{\E}X}{n-1}X+\frac{(b+4c)
div_{\E}JX}{n-1}JX,\\
b\n_{Z_{\lb}}\eta(Z_{\mu})=0, c\n_{Z_{\lb}}\eta(Z_{\mu})=0,\\
b(\n_{Z_{\lb}}\eta(Z_{\bar\mu})-\frac{
div_{\E}X}{2(n-1)}g_{\lb\bar\mu}+\frac{
div_{\E}JX}{2(n-1)}\0_{\lb\bar\mu})=0\\
c(\n_{Z_{\lb}}\eta(Z_{\bar\mu})-\frac{
div_{\E}X}{2(n-1)}g_{\lb\bar\mu}+\frac{
div_{\E}JX}{2(n-1)}\0_{\lb\bar\mu})=0\\
b\e(Z_{\lb})=0,\  b\e^*(Z_{\lb})=0\\
c(\e(Z_{\lb})+  \e^*(Z_{\lb}))=dc(Z_{\lb}).
\endgather$$

In particular in $U=\{x\in M:\kappa(x)\ne 0\}$ we obtain

$$\gather
\n a=\frac b{2(n-1)}\kappa \xi\tag 2.6\\
\n b=\frac {(b+4c)}{n-1}\kappa\xi,\\
\endgather$$
where $\xi$ is the principal section of the bundle $\De$. In the
rest of the paper we shall assume that the set $U=\{x:\kappa(x)\ne
0\}$ is non-empty and the set $B=\{x\in U:b(x)=0\}$ has an empty
interior. Thus we obtain
$$\gather \n_{Z_{\lb}}\eta(Z_{\mu})=0, \tag 2.7\\
\n_{Z_{\lb}}\eta(Z_{\bar\mu})-\frac{
div_{\E}X}{2(n-1)}g_{\lb\bar\mu}+\frac{
div_{\E}JX}{2(n-1)}\0_{\lb\bar\mu}=0\\
\e=0,\  \e^*=0.\endgather$$ It follows that the distribution
$\De_{U}$ is integrable. In $U$ we get
$$\n_{Z_{\lb}}\eta(Z_{\bar\mu})=\frac{\kappa}
{2(n-1)}g_{\lb\bar\mu}.\tag 2.8$$ Let $\xi$ be the principal
section of $\De$ defined in $U$ and let $\th=\xi^{\flat}$ be the
dual 1-form of $\xi$.  Then in $U$
$$da=\frac b{2(n-1)}\kappa \th,db=\frac {(b+4c)}{n-1}\kappa\th.\tag 2.9$$

It follows that the distribution $\Delta=\{X\in TU:\th(X)=0\}$
defined in $U$ is integrable. From (2.7) and (2.8) it follows that
$\n\th(JX,JY)=\n\th(X,Y)$ for $X,Y\in\G(\E)$. It is also clear
that
$$dJ\th_{|\E}=\frac{\kappa}{n-1}\0_{|\E}.$$
Thus the distribution $\E_{|U}$ is the so called $B_0$
-distribution defined in [G-M-1]. Consequently we obtain (see
[G-M-1], Lemma 5.1):
$$\n_X\th(Y)=\frac{\kappa}{2(n-1)}m(X,Y)+p\th(x)J\th(Y)-p^*J\th(X)J\th(Y),\tag
2.10$$ If dim $M=2n\ge 6$ we also have  $d\th=p \ \th\w J\th$,
$p=g(\n_{\xi}\xi,J\xi)=0$ and consequently (see Lemma 5.2 in
[G-M-1] and its proof):  $\n_{\xi}\xi=0, d\th=0$,
$$d\ln \kappa=-(\frac{\kappa}{n-1}+p^*)\th, dp^*\w\th=0.\tag
2.11$$ Thus
$$\n_X\th(Y)=\frac{\kappa}{2(n-1)}m(X,Y)-p^*J\th(X)J\th(Y),\tag
2.12$$ If $\m$ is a QCH K\"ahler manifold then one can show that
the Ricci tensor $\rho$ of $\m$ satisfies the equation
$$\rho(X,Y)=\lb m(X,Y)+\mu h(X,Y)\tag 2.13$$
where $\lb=\frac{n+1}2a+\frac  b4,\mu=\frac{n+1}2a+\frac{n+3}4b+c$
are eigenvalues of $\rho$ (see [G-M-1], Corollary 2.1 and Remark
2.1.) In particular the distributions $\E,\De$ are
eigendistributions of the tensor $\rho$ corresponding to the
eigenvalues $\lb,\mu$ of $\rho$.
\bigskip
{\bf 3.   Local holomorphic Killing vector field on $U$.} Let $\m$
be a QCH K\"ahler  manifold of dimension $2n\ge 6$.  We shall show
in this section that for every $x\in U$ there exists an open
neighborhood $V\subset U$ of $x$ and a function $f\in
C^{\infty}(V)$ such that $X_V=fJ\xi$ is a Killing vector field in
$V$, which we shall call a local special Killing vector field. Let
$V$ be a geodesically convex neighborhood of $x$ in $U$. Then $V$
is contractible. Note that the form $\phi=-p^*\th$, where
$p^*=g(\n_{J\xi}J\xi,\xi)$ is closed in $U$, since by (2.11),
$d\phi=-dp^*\w\th=0$. Consequently there exists a function $F\in
C^{\infty}(V)$ such that
$$dF=\phi=-p^*\th.\tag 3.1$$
Let $f=\exp\circ F$. From (2.12) it follows that
$$\n_XJ\th(Y)=\frac{\kappa}{2(n-1)}\0_m(X,Y)+p^*J\th(X)\th(Y),\tag
3.2$$

Now let $\psi=(fJ\xi)^{\flat}=fJ\th$. We shall show that
$\n_X\psi(Y)=-\n_Y\psi(X)$ which means that the field $fJ\xi$ is a
Killing vector field in $V$.We get

$$\n_X(fJ\th)(Y)=XfJ\th(Y)+f\frac{\kappa}{2(n-1)}\0_m(X,Y)+fp^*J\th(X)\th(Y),\tag
3.3$$ Since $Xf=fXF=-fp^*\th(X)$ we obtain
$$\gather \n_X(fJ\th)(Y)=f\frac{\kappa}{2(n-1)}\0_m(X,Y)-fp^*\th\w J\th(X,Y)=\tag
3.4\\f\frac{\kappa}{2(n-1)}\0_m(X,Y)-fp^*\om(X,Y),\endgather$$
which proves our claim. Note that if $F_1$ is another solution of
$(3.1)$ then $F_1=F+D$ for a certain constant $D\in\Bbb R$.  It
follows that $f_1=\exp F_1=CF$, where $C=\exp D$. Consequently
$X_1=f_1J\xi=CfJ\xi=CX$. Recall here the well known general fact,
that if $X,Y\in\frak{iso}(M)$ are Killing vector fields on
connected  Riemannian manifold $M$, $X\ne 0$ and $Y=fX$ for a
certain $f\in C^{\infty}(M)$ then $f$ is constant.

Let $\phi=f\xi^{\flat}$. Now it is clear that
$$\n_X\phi(Y)=\n_X(f\th)(Y)=f\frac{\kappa}{2(n-1)}m(X,Y)-fp^*h(X,Y)\tag 3.5$$
and $\n_X(f\th)(Y)=\n_Y(f\th)(X)$. Consequently $d(f\th)=0$.  It
follows that there exists a function $\tau\in C^{\infty}(V)$ such
that $f\xi=\n \tau$. Consequently $X_V=fJ\xi=J(\n \tau)$, which
means that $X_V$ is a holomorphic Killing vector field with a
K\"ahler potential $\tau$.
\medskip
{\bf 4. Special Jacobi fields along geodesics in $U$.}  Let
$c:[0,l]\rightarrow M$  be a unit geodesic such that
$c([0,l))\subset U$ and $c(l)\in K=\{x\in M:\kappa(x)=0\}$. A
vector field $C$, which is a Jacobi field along $c$ i.e. $\n_{\dot
c}^2C-R(\dot c,C)\dot c=0$, will be called a special Jacobi field
if there exists an open, geodesically convex neighborhood $V$ of
$c(0)$ such that $C(0)=X_V(c(0)),\n_{\dot c}C(0)=\n_{\dot
c}X_V(c(0))$. If $im\ c\cap V=c([0,\e))$ then it follows that
$X_V(c(t))=C(t)$ for all $t\in[0,\e)$. We have the following
lemma:
\medskip
{\bf Lemma 4.1.} {\it Let us assume that a vector  field $C$ along
a geodesic $c$ is a special Jacobi field along $c$. Then}

$$\lim_{t\rightarrow l}|C(t)|=0,$$
{\it and} $$g(\dot c,C)=0.$$

{\it Proof.} Let us note that $g(\n_{\dot c}C,\dot c)=0$ since
this property is valid for Killing vector fields. It follows that
the function  $g(\dot c,C)$ is constant. Let $k\in (0,l)$. Then
$c([0,k])\subset U$. For every  $t\in [0,k)$ there exists a
geodesically convex open neighborhood $V_t$ of the point $c(t)$
and a special Killing vector field $X_{V_t}=f_tJ\xi$ on $V_t$
defined in Section 3. The field $X_{V_t}$ is defined uniquely up
to a constant factor. From the cover $\{V_t\}:t\in [0,k]$ of the
compact set $c([0,k])$ we can choose a finite subcover $\{
V_{t_1},V_{t_2},...,V_{t_m}\}$. Let $c_i$ be the part of geodesic
$c$ contained in $V_i=V_{t_i}$, i.e. $im\ c\cap V_i=im\
c_i=c((t_i,t_{i+1}))$. We define the Killing vector field $X_i$ in
every $V_i=V_{t_i}$ by induction in such a way that $X_1=X_V$ on
$V_1\cap V$ and $X_i=X_{i+1}$ on $V_i\cap V_{i+1}$. Let
$X_i=f_iJ\xi$. Note that $C(t)=X_i\circ c(t)=f_iJ\xi\circ c(t)$
for $t\in (t_i,t_{i+1})$. Consequently, on $V_i$, $|C|=f_i$.  From
(2.11) and (3.1) it follows that
$$d\ln\kappa=d\ln f_i-\frac{\kappa}{n-1} \th.$$
Hence
$$\frac d{dt}\ln\kappa\circ c(t)=\frac d{dt}\ln |C(t)|-\frac{\kappa}{n-1}  \th(\dot c(t)),\tag 4.1$$
and

$$\frac d{dt}\ln\frac{\kappa\circ c(t)}{ |C(t)|}=-\frac{\kappa}{n-1}  \th(\dot c(t)).\tag 4.2$$
Consequently
$$\ln\frac{\kappa\circ c(k)}{ |C(k)|}-\ln\frac{\kappa\circ c(0)}{ |C(0)|}=-\frac1{n-1}\int_0^k\kappa \th(\dot
c(t))dt.\tag 4.3$$ Hence
$$\ln|C(k)|=\ln\kappa\circ c(k)-\ln\frac{\kappa\circ c(0)}{ |C(0)|}+\frac1{n-1}\int_0^k\kappa \th(\dot
c(t))dt.\tag 4.4$$ Note that
$$|\int_0^k\kappa \th(\dot
c(t))dt|\le \int_0^k|\kappa \th(\dot c(t))|dt\le \int_0^k\kappa
|\dot c(t)|dt\le \int_0^k\kappa dt.\tag 4.5$$ Let
$\kappa_0=sup\{\kappa(x):x\in c([0,l])\}$. From (4.4) it follows
that
$$\ln|C(k)|\le\ln\kappa\circ c(k)-\ln\frac{\kappa\circ c(0)}{ |C(0)|}+\frac1{n-1}\kappa_0l.\tag 4.6$$
Consequently
$$\limsup_{k\rightarrow l-}\ln|C(k)|\le\lim_{k\rightarrow l-}\ln\kappa\circ
c(k)-\ln\frac{\kappa\circ c(0)}{
|C(0)|}+\frac1{n-1}\kappa_0l=-\infty.\tag 4.7$$ From (4.7) it is
clear that $\lim_{k\rightarrow l-}|C(k)|=0.$ Since $|g(\dot
c(t),C(t))\le |C(t)|$ and $g(\dot c,C)$ is constant it follows
that $|g(\dot c(t),C(t))|\le \lim_{t\rightarrow l}|C(t)|=0$ which
means that $g(\dot c,C)=0$.$\k$
\medskip
{\bf 5. Global holomorphic Killing  vector field on $M$.} From now
on  we  assume that $\m$ is a complete QCH K\"ahler manifold,
$\dim M\ge 6$ and the set $U=\{x\in M:\kappa(x)\ne 0\}$ is
non-empty and $B$ has an empty interior.  Let $K=\{x\in
M:\kappa(x)=0\}$.
\medskip
{\bf Theorem 5.1.} {\it The set $U$ is connected and the set $K$
has an empty interior.}
\medskip
{\it Proof.} Let $U_1$ be a non-empty component of the set
$U=M-K$. Let $x_0\in U_1$ and let us assume the set $ int\ K\cup
(U-U_1)$ is non-empty. Let $x_1\in int\ K\cup(U-U_1)$. Then
$x_1=\exp_{x_0}lX$ for a certain unit vector $X\in T_{x_0}M$ and
$l>0$. Let $V\subset U_1$ be a geodesically convex neighborhood of
$x_0$ and $X_V=f_VJ\xi$ the local Killing vector field on $V$. Let
$C$ be the Jacobi vector field along the geodesic
$c(t)=\exp_{x_0}tX$ satisfying the initial
conditions:$C(0)=X_V(x_0),\n_{\dot c}C(0)=\n_{\dot c}X_V(x_0)$. It
follows that the field $C$ is a special Jacobi field along $c$. In
particular $g(X,X_V(x_0))=0$ since the geodesic $c$ meets $K$.
Since the set $int\ K\cup(U-U_1)$ is open there exists an open
neighborhood $W\subset int\ K\cup(U-U_1)$ of the point $x_1$. The
mapping $T_{x_0}M\ni Y\rightarrow \exp_{x_0}lY\in M$ is continuous
hence there exists an open neighborhood $P$ of $X$ in $T_{x_0}M$
such that $\exp_{x_0}lP\subset int\ K\cup(U-U_1)$. We can find a
vector $X_1\in P$ such that $g(X_1,X_V(x_0))\ne 0$. The field
$C_1$ along the geodesic $d(t)=\exp_{x_0}tX_1$ defined by the
initial conditions $C_1(0)=X_{V}(x_0),\n_{\dot d}C_1(0)=\n_{\dot
d(0)}X_V$ is a special Jacobi field along a geodesic $d$ which
meets $K$. It follows that $g(\dot d,C_1)=0$ along $d$. In
particular for $t=0$ we obtain $g(X_1,X_V(x_0))=0$ which is a
contradiction. Consequently $int\ K\cup(U-U_1)=\emptyset$.$\k$
\medskip
{\bf Corollary 5.2.} Let us  assume that $\m$ is a complete QCH
K\"ahler manifold, $\dim M\ge 6$,  the set $U=\{x\in
M:\kappa(x)\ne 0\}$ is non-empty and $B$ has an empty interior.
Then the distribution $\De$ is totally geodesic, i.e. $\n_XY\in
\G(\De)$ for every $X,Y\in\G(\De)$.
\medskip
{\it Proof.} Let $\n_XY=p_{\De}(\n_XY)+\a(X,Y)$ where
$\a(X,Y)=p_{\E}(\n_XY)$ and $\a$ is a section of the bundle
$\De^*\otimes\De^*\otimes\E$. Since $\De$ is integrable $\a$ is a
symmetric two-form.  We shall show that $\a=0$. The bundle
$\De_{|U}$ is spanned by the sections $\xi,J\xi$ and hence locally
by $-f\xi=JX,X=fJ\xi$ where $X$ is a special Killing field defined
in $dom f=V\subset U$. Note that $g(X,X)=f^2$ and
$\n_XX=-\frac12\n [g(X,X)]=-f\n f$. Thus (see (3.1))
$\n_XX=f^2p^*\xi\in \G(\De)$.
 Since $X$ is a holomorphic
Killing vector field it follows that $[X,JX]=L_X(JX)=J[X,X]=0$.
Hence $\n_XJX=\n_{JX}X$ and
$\n_XJX=J\n_XX=Xf\xi+f^2\n_{\xi}\xi=Xf\xi\in \G(\De)$ since
$\n_{\xi}\xi=0$. Finally $\n_{JX}JX=J(\n_{JX}X)=XfJ\xi\in
\G(\De)$. It follows that $\a(X,X)=\a(X,JX)=\a(JX,JX)=0$ and
$\a_{|U}=0$. Since $U$ is dense in $M$ it follows that $\a=0$ in
$M$. $\k$

\medskip
Let $x_0\in U$ and let $V=\exp_{x_0}W$ be a geodesically convex
open neighborhood of $x_0$, where $W\subset T_{x_0}M$ is a star
shaped open neighborhood of $0\in T_{x_0}M$. For every $X\in W$
let  $l(X)=\sup\{t:tX\in W\}$. Hence if the sphere $S_{\e}=\{X\in
T_{x_0}M:|X|=\e\}$ is contained in $W-\exp_{x_0}^{-1}(K)$ then
$W=\{tX:X\in S_{\e},t\in[0,l(X))\}$. In an open neighborhood
$V'\subset V-K$ of $x_0$ there is defined a special Killing vector
field $X_{V'}=f_{V'}J\xi$. Let $Z=X_{V'}(x_0)\in T_{x_0}M$ so that
$H=\{X\in T_{x_0}M:g_{x_0}(X,Z)=0\}$ is a hyperplane in
$T_{x_0}M$. Note that $\exp_{x_0}^{-1}(V\cap K)\subset H$. Let a
function $k:S_{\e}\rightarrow \Bbb R$ be defined as follows:
$k(X)=l(X)$ if $x\notin H$ and $k(X)=\inf\{t>0:\exp(tX)\in K\}$ if
$X\in H$. The set $W'=\{tX:X\in S_{\e},t\in[0,k(X))\}$ is open and
star shaped. Note that $V''=\exp_{x_0}W'\subset V-K\subset U$,
$V''$ is dense in $V$,  and $V''$ is contractible. It follows that
in $V''$ there is defined a special local  Killing vector field
$X_{V''}$. We can assume that $X_{V'}=X_{V''}$ on $V'$. Now we can
prove:

{\bf Lemma 5.3} {\it On every geodesically convex open set $V$ in
$M$ can be defined a  holomorphic  Killing vector field $X$ such
that for every open geodesically convex set $W\subset V\cap U$ the
restriction $X_{|W}$ is a special Killing vector field on $W$.}
\medskip
{\it Proof.} We shall use the notation introduced above. Since
$V''$ is contractible there exists a special Killing vector field
$X_{V''}$ defined on $V''$. Let us define a differentiable field
$X$ on $V$ by the formula: $X(\exp_{x_0}u)=J_u(1)$ where $u\in W$
and $J_u$ is a Jacobi vector field along a geodesic
$c(t)=\exp_{x_0}(tu)$ satisfying the initial
conditions:$J_u(0)=X_{V''}(x_0),\n_{\dot c}J_u(0)=\n_{u}X_{V''}$.
It is clear that $X$ is a differentiable vector field and that
$X_{|V''}=X_{V''}$. Since the set $V''$ is dense in $V$ it follows
that $X$ is a Killing vector field in $V$. In fact if we write
$T=\n X$ then for every $Y,Z\in\frak X(V)$ we have
$g(TY,Z)=-g(Y,TZ)$ on $V''$ and both sides are differentiable
functions on $V$. Thus the relation remains valid on $V$ which
means that $X$ is a Killing vector field on $V$. Note that $X=0$
on $V\cap K$.  Since  equation (3.5) is valid on the open, dense
subset $V-K$ of $V$ it follows that the form $\phi=-(JX)^{\flat}$
satisfies the relation $\n_Y\phi(Z)=\n_Z\phi(Y)$ for every
$Y,Z\in\frak X(V)$. Consequently $JX= -\n \tau$ for a certain
function $\tau\in C^{\infty}(V)$ and $X=J\n \tau$ on $V$. Since
$\n_{JY}\phi(JZ)=\n_Y\phi(Z)$ on $V''$ and hence on $V$ for every
$Y,Z\in \frak X(V)$ it is clear that $X$ is holomorphic. $\k$
\medskip
Let $V_1,V_2$ be two open, geodesically convex sets. Then the set
$V_1\cap V_2-K$ is connected. The proof of this is similar to the
proof of Theorem 5.1. Since the sets $V_1''\subset V_1-K$,
$V_2''\subset V_2-K$ are contractible there exist special Killing
vector fields $X_1=f_1J\xi, X_2=f_2J\xi$ on these sets. These
Killing fields can be extended on $V_1-K,V_2-K$ in such a way that
$f_i=\exp F_i$ where $F_i$ satisfy   equation (3.1). Consequently
$d(F_1-F_2)=0$ on the connected set $V_1\cap V_2-K$. Hence
$F_1=F_2+D$ in $V_1\cap V_2-K$ for a constant $D\in \Bbb R$.
Consequently $X_1=C_{12}X_2$ where $C_{12}=\exp D$. The fields
$X_1,X_2$ can be extended to the Killing fields on $V_1,V_2$
respectively such that $X_{i|K\cap V_i}=0$. It is clear that for
these extensions which we also denote by $X_1,X_2$ the equation
$$X_1=C_{12}X_2$$
holds on $V_1\cap V_2$.
\bigskip
Let us recall the definition of a special K\"ahler-Ricci potential
([D-M-1],[D-M-2]).
 \medskip {\bf Definition.} {\it A nonconstant
function $\tau\in C^{\infty}(M)$, where $\m$ is a K\"ahler
manifold, is called a special K\"ahler-Ricci potential if the
field $X= J(\n \tau)$ is a Killing vector field and, at every
point with $d\tau\ne 0$ all nonzero tangent vectors orthogonal to
the fields $X,JX$ are  eigenvectors of both $\n d\tau$ and the
Ricci tensor $\rho$ of $\m$.}
\medskip
Now we shall prove:
\bigskip
{\bf Theorem 5.4.  }{\it   Let $\m$ be a complete  QCH K\"ahler
manifold of dimension $2n\ge 6$. Let $int B=\emptyset$ and $U\ne
\emptyset$. If $H^1(M,\Bbb R)=0$ then there exists on $M$ a
non-zero holomorphic Killing vector field $X=J(\n \tau)$ with a
special K\"ahler-Ricci potential $\tau$.}
\medskip
{\it Proof.} Let $\{V_i\}:i\in I$ be a cover of $M$ by
geodesically convex, open sets $V_i$. Let $X_i$ be a Killing
vector field on $V_i$ constructed in Lemma 5.2. Let $V_i\cap
V_j\ne\emptyset$. Then there exist constants $C_{ij}>0$ such that
$X_i=C_{ij}X_j$ on $V_i\cap V_j$. These constants satisfy the
co-cycle condition $C_{ij}C_{jk}C_{ki}=1$. Consequently the
constants $D_{ij}=\ln C_{ij}\in \Bbb R$ satisfy the co-cycle
condition $D_{ij}+D_{jk}+D_{ki}=0$. It follows that $\{D_{ij}\}$
is a co-cycle in the first \v{C}ech cohomology group
$\breve{H}^1(\{V_i\},\Bbb R)$. Since $\{V_i\}$ is a good cover of
$M$ it follows that $\breve{H}^1(\{V_i\},\Bbb
R)=\breve{H}^1(M,\Bbb R)=0$. Consequently there exists a co-cycle
$\{D_i\}\in Z^0(\{V_i\},\Bbb R)$ such that
$\{D_{ij}\}=\dl(\{D_i\})$. This means that $D_{ij}=D_j-D_i$. Let
$C_i=\exp D_i$. Then $C_{ij}=\frac{C_j}{C_i}$. Let us define the
field $X$ on $M$ by the formula $X_{|V_i}=C_iX_i$. Then it is
clear that $X$ is a well defined, global vector field and
$X\in\frak X(M)$. Since $X_{|V_i}$ is a Killing vector field on
every $V_i$ it follows that $X$ is a Killing vector field. Now let
$\phi=-(JX)^{\flat}$. Then $d\phi=0$, since this equation is
satisfied on every $V_i$. On the other hand the first de Rham
group of $M$ vanishes: $H^1(M,\Bbb R)=\breve{H}^1(M,\Bbb R)=0$.
Consequently there exists a function $\tau\in C^{\infty}(M)$, such
that $\phi=d\tau$. Note also that $\n d\phi$ is Hermitian, which
means that $X$ is  holomorphic.   Thus $X=J(\n \tau)$ is a
holomorphic Killing vector field with a Killing potential $\tau$.
Note that in view of (2.12) and (3.5) the special Killing field
constructed by us is a Killing vector field with a special
K\"ahler-Ricci potential $\tau$. $\k$
\medskip
\bigskip
{\bf Corollary  5.5.  }{\it   Let $\m$ be a complete  QCH K\"ahler
manifold of dimension $2n\ge 6$. Let $int B=\emptyset$,  $U\ne
\emptyset$ and let $(\tilde M,\tilde g)$ be the Riemannian
universal covering space of $\m$.  Then there exists on $(\tilde
M,\tilde g)$ a non-zero holomorphic Killing vector field with a
special K\"ahler-Ricci  potential.}

\bigskip
{\bf 6. Construction of QCH K\"ahler manifolds.} In our
construction we shall follow B\'erard Bergery (see [Ber], [S])
rather then Derdzi\'nski and Maschler, although we shall use the
classification theorem by Derdzi\'nski and Maschler ([D-M-1]) to
classify QCH K\"ahler manifolds. Note that these two approaches
are equivalent (see [D-M-2]).  Let $(N,h,J)$ be a  simply
connected  K\"ahler Einstein manifold, which is not Ricci flat and
$\dim N=2m>2$. Let $s\ge0,L>0,s\in \Bbb Q,L\in \Bbb R$ and
$r:[0,L]\rightarrow \Bbb R$ be a positive, smooth function on
$[0,L]$ with $r'(t)>0$ for $t\in (0,L)$, which is even at $0$ and
$L$, i.e. there exists an $\e>0$ and  even, smooth functions
$r_1,r_2:(-\e,\e)\rightarrow \Bbb R$ such that $r(t)=r_1(t)$ for
$t\in[0,\e)$ and $r(t)=r_2(L-t)$ for $t\in(L-\e,L]$. If $s\ne0$
then it is clear that the function $f=\frac2srr'$ is positive on
$(0,L)$ and $f(0)=f(L)=0$. Let $P$ be a circle bundle over $N$
classified by the integral cohomology class $\frac s2c_1(N)\in
H^2(N,\Bbb R)$ if $c_1(M)\ne 0$. Let $q$ be the unique positive
integer such that $c_1(N)=q\a$ where $\a\in H^2(N,\Bbb R)$ is an
indivisible integral class. Then
$$s=\frac{2k}q; k\in\Bbb Z.$$ It is known that $q=n$ if
$N=\Bbb{CP}^{n-1}$ (see [Bes], p.273). Note that
$c_1(N)=\{\frac1{2\pi}\rho_N\}=\{\frac{\tau_N}{4m\pi}\0_N\}$ where
$\rho_N=\frac{\tau_N}{2m}\0_N$ is the Ricci form of $(N,h,J)$,
$\tau_N$ is the scalar curvature of $(N,h)$ and $\0_N$ is the
K\"ahler form of $(N,h,J)$. We can assume that $\tau_N=\pm 4m$. In
the case $c_1(N)=0$ we shall assume that $(N,h,J)$ is a Hodge
manifold, i.e. the cohomology class $\{\frac s{2\pi}\0_N\}$ is an
integral class. On the bundle $p:P\rightarrow N$ there exists a
connection form $\th$ such that $d\th=sp^*\0_N$ where
$p:P\rightarrow N$ is the bundle projection. Let us consider the
manifold $(0,L)\times P$ with the metric
$$g=dt^2+f(t)^2\th^2+r(t)^2p^*h,\tag 6.1$$ if $s\ne0$ and the
metric
$$g=dt^2+f(t)^2\th^2+p^*h,\tag 6.2$$
if $s=0$. The metric (6.2) is a K\"ahler product metric.
 The metric $6.1$ is
K\"ahler if and only if $f=\frac{2rr'}s$. We shall prove it in
section 8. It is known that the metric (6.1) extends to a metric
on a sphere bundle $M=P\times_{S^1}\Bbb{CP}^1$ if and only if the
function $r$ is positive and smooth on $(0,L)$, even at the points
$0,L$, the function $f$ is positive, smooth an odd at the points
$0,L$ and additionally $$f'(0)=1,f'(L)=-1.\tag 6.3$$ If
$f=\frac{2rr'}s$ for $r$ as above, then (6.3) means that
$$2r(0)r''(0)=s,2r(L)r''(L)=-s.\tag 6.4$$

{\bf 7.  Circle bundles.  } Let $(N,h,J)$ be a K\"ahler manifold
with integral class  $\{\frac s{2\pi}\0_N\}$, where $s\in\Bbb Q$
and let $p:P\rightarrow N$ be a circle bundle with a connection
form $\th$ such that $d\th=s\0_N$ (see [K]). We shall assume that
$(N,h)$ is Einstein if $dim\ N>2$.  Let us consider a Riemannian
metric $g$ on $P$ given by
$$g=a^2\th\otimes\th+b^2p^*h\tag 7.1$$
where $a,b\in \Bbb R$. Let $\xi$ be the fundamental vector field
of the action of $S^1$ on $P$ i.e. $\th(\xi)=1, L_{\xi}g=0$. It
follows that $\xi\in\frak{iso}(P)$ and $a^2\th=g(\xi,.)$.
Consequently
$$a^2d\th(X,Y)=2g(TX,Y)\tag 7.2$$
for every $X,Y\in\frak X(P)$ where $TX=\n_X\xi$. Note that
$g(\xi,\xi)=a^2$ is constant, hence $T\xi=0$. On the other hand
$d\th(X,Y)=sp^*\0_N(X,Y)=sh(Jp(X),p(Y))$.  Note that there exists
a tensor field $\tilde J$ on $P$ such that $\tilde J\xi=0$ and
$\tilde J(X)=(JX_*)^H$ where $X=X_*^H\in TP$ is the horizontal
lift of $X_*\in TN$ (i.e. $\th(X_*^H)=0$) and $X_*=p(X)$. Indeed
$L_{\xi}T=0$ and $T\xi=0$ hence $T$ is the horizontal lift of the
tensor $\tilde T$. Now $\tilde J=\frac{2b^2}{sa^2}\tilde T$. Since
$T\xi=0$ we get $\n T(X,\xi)+T^2X=0$ and $R(X,\xi)\xi=-T^2X$. Thus
$g(R(X,\xi)\xi,X)=||TX||^2$ and
$$\rho(\xi,\xi)=||T||^2=\frac{sa^4}{4b^4}2m.$$ Consequently
$$\lb=\rho(\frac{\xi}a,\frac{\xi}a)=\frac1{a^2}||T||^2=\frac{s^2a^2}{4b^4}m.\tag 7.3$$

We shall compute the O'Neill tensor $A$ (see [ON]) of the
Riemannian submersion $p:(P,g)\rightarrow (N,b^2h)$. We have
$$A_EF=\Cal V(\n_{\Cal H E}\Cal H F)+\Cal H(\n_{\Cal H E}\Cal V
F).$$

Let us write $u=\Cal V(\n_{\Cal H E}\Cal H F)$ and $v=\Cal
H(\n_{\Cal H E}\Cal V F)$.  The vertical component of a field $E$
equals $\th(E)\xi$. If $X,Y\in\Cal H$ then
$$g(\n_XY,\frac1a\xi)=\frac1a
(Xg(Y,\xi)-g(Y,\n_X\xi))=-\frac1ag(TX,Y)=\frac1ag(X,TY).\tag 7.4$$
Hence
$u=\frac1ag(E-\th(E)\xi,T(F-\th(F)\xi)\frac{\xi}a=\frac1{a^2}g(E,TF)\xi.$
Note that $\Cal H(\n_Xf\xi)=f\Cal H(\n_X\xi)=fTX$, hence $$v=\Cal
H(\n_{\Cal H E}\Cal V F)=\th(F)T(E)=\frac1{a^2}g(\xi,F)TE.$$
Consequently
$$A_EF=\frac1{a^2}(g(E,TF)\xi+g(\xi,F)TE).\tag 7.5$$
If $U,V\in \Cal H$ then
$$||A_UV||^2=\frac1{a^2}g(E,TF)^2=\frac{s^2a^2}{4b^4}g(E,\tilde J
F)^2.$$ If $E$ is horizontal and $F$ is vertical then
$$A_EF=\frac1{a^2}g(\xi,F)TE.\tag 7.6$$
Hence $A_E\xi=TE$ and $||A_E\xi||^2=||TE||^2=\frac{s^2a^4}{4b^4}$.
It follows that
$$K(P_{E\xi})=\frac{s^2a^2}{4b^4},$$
where $K(P_{EF})$ denotes the sectional curvature of the plane
generated by vectors $E,F$.  If $E,F\in\Cal H$ then
$$K(P_{EF})=K_*(P_{E_*F_*})-\frac{3g(E,TF)^2}{a^2||E\w F||^2},$$
where $E_*$ denotes the projection of $E$ on $M$ i.e. $E_*=p(E)$.
Thus
$$K(P_{EF})=\frac1{b^2}K_0(P_{E_*F_*})-\frac{3s^2a^2g(E,\tilde JF)^2}{4b^4||E\w F||^2},\tag 7.7$$
where $K_0$ stands for the sectional curvature of the metric $h$
on $N$. Applying this we get for any $E\in\Cal H$ the formula for
the Ricci tensor $\rho$ of $(M,g)$:
$$\rho(E,E)=\frac1{b^2}\rho_0(bE_*,bE_*)-\frac{3s^2a^2}{4b^4}+\frac{s^2a^2}{4b^4},$$
where $\rho_0$ is a Ricci tensor of $(M,h)$. Hence
$$\mu=\frac{\mu_0}{b^2}-\frac{s^2a^2}{2b^4},$$
where $\rho_0=\mu_0g_0$. Now we shall find a formula for
$R(X,\xi)Y$ where $X,Y\in \Cal H$. We have $R(X,\xi)Y=\n T(X,Y)$
and $$\gather \n T(X,Y)=\n_X(T(Y))-T(\n_XY)=\n^*_{X_*}(\tilde
TY^*)+\frac12\Cal V[X,TY]\tag 7.8\\-(\tilde
T(\n^*_{X_*}Y_*))^*=\frac12\Cal
V[X,TY]=-\frac12sp^*\0_N(X,TY)\xi=-\frac{s^2a^2}{4b^2}h(X_*,Y_*)\xi\endgather$$
Consequently $R(X,Y,Z,\xi)=0$ for $X,Y,Z\in\Cal H$,and
$$R(X,\xi,Y,\xi)=-\frac{s^2a^4}{4b^2}h(X_*,Y_*).\tag 7.9$$

{\bf 8.   Riemannian submersion $  p:(0,L)\times P\rightarrow
(0,L)$.}  In this case the O'Neill tensor $A=0$. We shall compute
the O'Neill tensor $T$ (see [ON]). Let us denote by  $Y^*$ the
horizontal lift of the vector $Y\in TN$ with respect to the
Riemannian  submersion $p_N:P\rightarrow N$ i.e.
$p_N(Y^*)=Y,g(Y^*,\xi)=0$. Let $H=\frac d{dt}$ be the horizontal
vector field for this submersion and $\De$ be the distribution
spanned by the vector fields $H,\xi$. If $U,V\in \Cal V$ and
$g(U,V)=0$ then $T(U,V)=0$. Let  $U\in \Cal V,g(U,\xi)=0$ and
$U=U_*^*$ with $h(U_*,U_*)=1$, then the following formula holds
$$T(U,U)=-rr'H.\tag 8.1$$
In fact $2g(\n_UV,H)=-Hg(U,V)=-2rr'h(U_*,V_*)$ if $U=V$ or $0$ if
$g(U,V)=0$.
 We also have
$$T(\xi,\xi)=-ff'H.\tag 8.2$$
Now we shall prove that the almost complex structure defined by
$$JH=\frac1f\xi,JX=(J_*X_*)^*\text{   for  } X=(X_*)^*\in\E=\De^{\perp}$$
where $X_*\in TN$, is a K\"ahler structure with respect to the
metric $g$. We have for horizontal lifts $X,Y\in\frak
X(P)\subset\frak X((0, L)\times P)$ of the fields $X_*,Y_*\in\frak
X(N)$ ( with respect to the submersion described in the Section
7):
$$\gather \n
J(Y,X)=\n_Y(JX)-J(\n_YX)=\n^*_{Y_*}(J_*(X_*))^*-\frac12d\th(Y,JX)\xi\\
+T(Y,JX)-J(\n^*_{Y_*}(X_*)^*-\frac12d\th(Y,X)\xi+T(Y,X))=\\
-\frac12sh(JY,JX)fJH-rr'h(Y,JX)H-\frac12sh(JY,X)fH+h(X,Y)rr'JH=0
\endgather$$
if and only if $f=\frac{2rr'}s$.  Since the distribution $\De$ is
totally geodesic and two-dimensional  it is clear that $\n
J(X,Y)=0$ if $X,Y\in\G(\De)$.  Now we shall show that
$$\n J(H,X)=\n J(X,H)=0\text {  for   } X\in \E.$$
It is easy to show that $\n_XH=\n_HX=\frac{r'}rX$ and
$$\n_X(JH)=\n_X(\frac1f\xi)=\frac1fT(X)=\frac{sf}{2r^2}JX.$$
On the other hand
$$\n_X(JH)=\n J(X,H)+J(\n_XH)=\n J(X,H)+\frac {r'}rJ(X).$$
Thus $\n J(X,H)=0$ if $f=\frac{2rr'}s$.  Similarly
$$\n_H(JX)=\n_{JX}H=\frac{r'}rJX=\n J(H,X)+J(\n_HX)=\n
J(H,X)+\frac{r'}rJX$$ and $\n J(H,X)=0$. Note that the K\"ahler
form $\0=fdt\w\th+r^2p^*\0_N$ of almost Hermitian manifold
$((0,L)\times P,g,J)$ is closed, which means that the structure
$J$ is almost K\"ahler. Thus $\n J(JX,Y)=-J\n J(X,Y)$ and
consequently $\n_{\xi}J=0$ which finishes the proof. Let
$U,V,W\in\Cal V$ and $g(U,\xi)=g(V,\xi)=g(W,\xi)=0$. Then
$$R(U,V,\xi,W)=\hat
R(U,V,\xi,W)-g(T(U,\xi),T(V,W))+g(T(V,\xi),T(U,W))=0.$$ From
O'Neill formulae it follows also that
$$R(JH,U,V,JH)=0$$ if $g(U,V)=0$ and
$$R(JH,U,U,JH)=\frac{s^2f^2}{4r^4}+\frac{f'r'}{fr},\tag 8.3$$ for
a unit vector field $U$ as above. Note also that the distribution
$\De$ spanned by the vector fields $\xi,H$ is totally geodesic.
Consequently if $X,Y,Z\in\G(\De)$ and $V$ is as above then
$$R(X,Y,Z,V)=0.\tag 8.4$$

{\bf Theorem 8.1.} {\it The K\"ahler metric
$g=dt^2+(\frac{2rr'}s)^2\th^2+r^2p^*h$ on the manifold
$(0,L)\times P$ has quasi-constant holomorphic sectional curvature
if and only if dim$N=2$ or the K\"ahler manifold $(N,h)$ has
constant holomorphic sectional curvature and dim$N\ge 4$.}

\medskip
{\it Proof.} Let $X=X^*+\a H+\be JH$, where $X^*\in
\E=\De^{\perp}$. Let us assume that $g(X,X)=1$. Note that
$|X_{\De}|^2=\a^2+\be^2$ and $|X^*|^2=1-|X_{\De}|^2$ and
$JH=\frac1f\xi$ where $f=\frac{2rr'}s$. We also have

$$\gather R(X^*+\a H+\be JH,JX^*+\a JH-\be H,JX^*+\a JH-\be H,X^*+\a H+\be
JH)=\\a^4R(H,JH,JH,H)-a^2\be^2R(H,JH,H,JH)+\be^4R(JH,H,H,JH)\\
-\a^2\be^2R(JH,H,JH,H)+\a^2R(H,JX^*,JH,X^*)+\a^2R(H,JX^*,JX^*,H)\\+\a^2R(X^*,JH,JX^*,H)
+\a^2R(X^*,JH,JH,X^*) +\be^2R(JH,JX^*,-H,X^*)\\
\be^2R(JH,JX^*,JX^*,JH)+\be^2R(X^*,H,H,X^*)+\be^2R(X^*,-H,JX^*,JH)\\
+2\a^2R(H,JH,JX^*,X^*)-2\be^2R(JH,H,JX^*,X^*)\\
+R(X^*,JX^*,JX^*,X^*)=|X_{\De}|^4R(H,JH,JH,H)+8|X_{\De}|^2R(H,JX^*,JX^*,H)\\
+R(X^*,JX^*,JX^*,X^*)=|X_{\De}|^4R(H,JH,JH,H)+\\8|X_{\De}|^2(1-|X_{\De}|^2)R(JH,\tilde
X^*,\tilde X^*,JH) +(1-|X_{\De}|^2)^2R(\tilde X^*,J\tilde
X^*,J\tilde X^*,\tilde X^*)
\endgather$$
where $\tilde X^*=\frac1{|X^*|}X^*$ if $X^*\ne 0$ or $\tilde
X^*=0$ if $X^*=0$.

From (8.5) it follows that $R(X^*+\a H+\be JH,JX^*+\a JH-\be
H,JX^*+\a JH-\be H,X^*+\a H+\be JH)$ depends only on $|X_{\De}|$
and the point $x\in (0,L)\times P$ if and only if $R(\tilde
X^*,J\tilde X^*,J\tilde X^*,\tilde X^*)$ does not depend on the
unit vector $\tilde X^*$. It follows that $dim N=2$ or $dim N\ge
4$ and  at every point $x\in N$  the holomorphic sectional
curvature of $(N,h)$ is constant. In fact from O'Neill formulae
([ON]) it follows that $$\gather R(\tilde X^*,J\tilde X^*,J\tilde
X^*,\tilde X^*)=R_*(\tilde X^*,J\tilde X^*,J\tilde X^*,\tilde
X^*)-g(T(\tilde X^*,\tilde X^*),T(J\tilde X^*,J\tilde
X^*))\\+2g(A(\tilde X^*,J\tilde X^*),A(J\tilde X^*,\tilde
X^*))-g(A(J\tilde X^*,\tilde X^*),A(J\tilde X^*,\tilde
X^*))\\=R_*(\tilde X^*,J\tilde X^*,J\tilde X^*,\tilde
X^*)-4\frac{(r')^2}{r^2}=\frac{c_0}{r^2}-4\frac{(r')^2}{r^2},\endgather$$
where we used the formula $A(E,F)=\frac s{2r^2}g(E,\tilde JF)\xi$.
 Note that
if $(N,h,J)$ has constant holomorphic sectional curvature $c_0$
then the coefficient $a$ in the formula (2.1) for the tensor $R$
equals $a=\frac{c_0}{r^2}-4\frac{(r')^2}{r^2} $ and is
non-constant. Note that $div_{\E}H=2(n-1)\frac{r'}r$ and
$div_{\E}\xi=0$. In particular $\kappa=2(n-1)\frac{r'}r\ne 0$ on
an open and dense subset.
 $\k$
\medskip
{\it Remark.} Note that in the case $s=0$ we also get a K\"ahler
manifold with QCH metric and a special Killing Ricci potential.
However in this case $\kappa=0$ on the whole of $(0,L)\times P$.
Let $g_0$ be any smooth, Riemannian metric on $\Bbb{CP}^1$. If $h$
is a metric of constant holomorphic sectional curvature on
$\Bbb{CP}^n$, then the product
$(\Bbb{CP}^1\times\Bbb{CP}^{n-1},g_0 \times h)$ is a compact,
simply connected manifold QCH K\"ahler manifold (with respect to
the distribution  whose leaves are submanifols
$\Bbb{CP}^1\times\{x\}$, for $x\in \Bbb{CP}^{n-1}$), having
$\kappa=0$ and, in general, not admitting a special K\"ahler-Ricci
potential.
\medskip
{\bf 9. QCH K\"ahler  manifolds with real analytic Riemannian
metric.}  We start with :
\medskip
{\bf Lemma 9.1.} {\it Let us assume that $\m$ is a (connected,
complete) QCH K\"ahler  manifold with a real analytic metric $g$
and dim$\ M\ge6$. Then int$\ B=\emptyset$ or $b=0$ and $\m$ has a
constant holomorphic sectional curvature.}

 {\it Proof.} We can find local coordinates $(x_1,x_2,...,x_{2n})$ such that
 $g_{\a\be}$ are real analytic, where $g_{\a\be}=g(\p_{\a},\p_{\be})$ and
 $\p_{\a}=\frac{\p}{\p x_{\a}}$. Analogously
$\rho_{\a\be}=\rho(\p_{\a},\p_{\be})$ are real analytic functions.
If $g(SX,Y)=\rho(X,Y)$, where $\rho$ is the Ricci tensor of $\m$,
then $S^{\be}_{\a}=g^{\g\be}\rho_{\a\g}$ where
$S\p_{\a}=S_{\a}^{\be}\p_{\be}$. It follows that the functions
$S_{\a}^{\be}$ are analytic. Consequently the polynomial
$\phi(t)=det(S-tI)$ has real analytic coefficients. Thus if $b=0$
in open set $V\subset U$ then also $c=0$ on $V$ and $S=\lb I$ on
$V$  for a constant $\lb\in\Bbb R$.  This means that
$\phi(t)=(\lb-t)^{2n}$ in $V$ and consequently on the whole of
$M$. Thus $S=\lb I$ on the whole of $M$. It follows  (see (2.13))
 that $\frac{n+1}2a+\frac b4$ is constant on the whole of $M$ and
$\frac{n+2}4b+c=0$ on the whole of $M$. Consequently
$b+4c=-(n+1)b$. From (2.9) we get in $U-B$:

$$  \xi(\ln |b|)=-\frac{n+1}{n-1}\kappa         \tag 9.1$$
Let $t\rightarrow c(t)$ be a unit geodesic joining a point
$c(0)=x\in U-B$ with a point $c(l)=y\in B$. Note that $|\dot c(\ln
|b|)|\le |\xi(\ln |b|)|=\frac{n+1}{n-1}\kappa$. Thus we obtain
$\ln |b|\circ c(t)-\ln |b|\circ
c(0)\ge-\frac{n+1}{n-1}\int_0^t\kappa\circ c(s)ds\ge
-\frac{n+1}{n-1}\kappa_0l$ where $0\le \kappa\circ c(s)\le
\kappa_0=\sup_{c([0,l])}\kappa $. It is clear now that
$U-B=\emptyset$ and $b=0$ in $U$. It follows that in the case of
analytic metric the assumption $int\ B\ne \emptyset$ implies
$b=0,c=0$ on the whole of $U$ and hence on $M$ which means that
$\m$ has constant holomorphic sectional curvature.$\k$
\medskip
Now we shall prove
\medskip
{\bf Theorem 9.2. }{\it Let $\m$ be a compact, simply connected
QCH K\"ahler manifold of dimension $2n\ge 6$. If $\kappa\ne 0$ and
$int\ B=\emptyset$ then $M$ is a $\Bbb{CP}^1$-bundle over
$\Bbb{CP}^{n-1}$ with a metric homothetic to the metric (6.1)
where $s\ne 0,f=\frac{2rr'}s,s=\frac{2k}n,k\in\Bbb N$, $r$
satisfies the boundary conditions (6.4) and $h$ is a K\"ahler
metric of constant holomorphic sectional curvature on
$\Bbb{CP}^{n-1}$. }
\medskip
{\it Proof.} From Theorem 5.3 it follows that $\m$ is a compact
K\"ahler  manifold admitting a special Killing-Ricci potential.
From the classification of such manifolds given by Derdzi\'nski
and Maschler ( see  [D-M-1]) it follows that $M$ is a holomorphic
$\Bbb{CP}^1$-bundle over a compact, simply connected
Einstein-K\"ahler manifold $(N,h)$, i.e. $r$ is given by formulae
(6.1) or (6.2) with $f,r$ satisfying the initial conditions. Since
$\kappa\ne 0$ it follows that $s\ne 0$ and the metric is
homethetic to the metric (6.1) with $f=\frac{2rr'}s$ where $r$ is
a smooth function, positive on $(0,L)$, even at $0,L$ and with
$r'>0$ satisfying the boundary conditions (6.4). As $\m$ is a QCH
K\"ahler manifold it follows that $(N,h)$ is compact, simply
connected K\"ahler manifold of constant holomorphic curvature,
which means that $N=\Bbb{CP}^{n-1}$ with a standard Fubini-Study
metric.$\k$

\medskip
{\bf Theorem 9.3. }{\it Let $\m$ be a compact, simply connected
QCH K\"ahler manifold of dimension $2n\ge 6$  with an analytic
Riemannian metric $g$. If $\kappa\ne 0$  then $M$ is a
$\Bbb{CP}^1$-bundle over $\Bbb{CP}^{n-1}$ with a metric homothetic
to the metric (6.1) where $s\ne 0,f=\frac{2rr'}s,
s=\frac{2k}n,k\in\Bbb N$ with $r\in C^{\om}(\Bbb R)$, $r$
satisfies the boundary conditions (6.4) and $h$ is a K\"ahler
metric of constant holomorphic sectional curvature on
$\Bbb{CP}^{n-1}$ with scalar curvature $4(n-1)$ or $\m$ is a
projective complex space $\Bbb{CP}^n$ with a metric homothetic to
the standard Fubini-Study metric. }
\medskip
{\it Proof.} If the metric $g$ is analytic then  $int\
B=\emptyset$ or $\{x:b(x)=0,c(x)=0\}=M$. In the first case we
apply Theorem 9.2 and in the second case $\m$ has constant
holomorphic sectional curvature at every point $x\in M$ which
means that $\m$ has constant holomorphic sectional curvature (see
[K-N], Chapter 9, Theorem 7.5, p.158 of Russian translation).
Since $M$ is compact and simply connected it follows (see [K-N],
Chapter 9, Theorems 7.8, 7.9 pp.160-161) that $M=\Bbb{CP}^n$ with
a metric homothetic to the standard Fubini-Study metric. $\k$
\medskip
{\it  Remark. } Note, that the tangent bundle $T\Bbb{CP}^n$ may
not admit any complex line subbundle. This is the case if $n=2$
$(T\Bbb{CP}^2$ does not contain any two-dimensional, oriented real
subbundle, see [H-H]). Thus $\Bbb{CP}^2$ does not admit any
$J$-invariant two-dimensional global distribution $\De$ and cannot
be a QCH K\"ahler manifold.
\bigskip
Now we show that there are uncountably many  analytic functions
$r$ satisfying the boundary conditions (6.4). Note first that the
function $(r')^2$ is a function of $r$ alone, i.e. $(r')^2=P(r)$
for a certain function $P$ which is smooth if $r$ is smooth. We
shall give a family of analytic functions $P$ which are
polynomials of degree 3, parameterized by the real numbers $x,y$
where $im\ r=[x,y]$ which give rise to functions $r$ satisfying
the boundary conditions  (6.4) and thus giving examples of QCH
K\"ahler manifolds. In that way we show one of the  method of
constructing  such functions. Let $0<x<y$ and let us consider a
polynomial
$$P(t)=-\frac{s(x+y)}{y-x}+\frac{s(x^2+3sxy+y^2)}{xy(y-x)}t-\frac{2s(x+y)}{xy(y-x)}t^2+\frac
s{xy(y-x)}t^3.\tag 9.2$$ Then $P(x)=P(y)=0$ and $P(t)>0$ for $t\in
(x,y)$ and $$ P'(x)x=s,P'(y)y=-s.\tag 9.3$$  Let us consider the
ordinary differential equation
$$r''=\frac12P'(r)\tag 9.4$$ with the initial conditions $r(0)=x,r'(0)=0$.
Then any solution $r$ of $(9.4)$ is an analytic function $r\in
C^{\om}(\Bbb R)$. Let $L>0$ be the first point such that $r(L)=y$.
Then for $t\in (0,L)$ we have
$$r'(t)=\sqrt{P(r)}.\tag 9.5$$
Note that
$$L=\int_x^y\frac{dt}{\sqrt{P(t)}}.\tag 9.6$$
The function $r$ is even at both points $0,L$. In fact we have
$r'(0)=0,r'(L)=0$ and $r$ satisfies (9.4). It follows that
$r^{(2k-1)}(0)=r^{(2k-1)}(L)=0$ for $k\in\Bbb N$. Since $r\in
C^{\om}(\Bbb R)$ it follows that $r$ is even at $0,L$. In view of
$(9.5)$ $r'(t)>0$ for every $t\in (0,L)$. We also have
$$\gather 2r(0)r''(0)=2x\frac12P'(x)=xP'(x)=s\\
2r(L)r''(L)=2y\frac12P'(y)=yP'(y)=-s.\endgather$$

It follows that the metric $g=dt^2+(\frac{2rr'}s)^2\th^2+r^2p^*h$
on $(0,L)\times P$ where $L$ is given by (9.6) and $r$ is the
solution of (9.4) extends to an analytic metric on the
$\Bbb{CP}^1$-bundle $M=P\times_{S^1}\Bbb{CP}^1$. Note that the
special Killing-Ricci potential is $\tau=\frac{r^2}s$ and the
holomorphic Killing vector field corresponding to it is $\xi$ with
$g(\xi,\xi)=f^2=(\frac{2rr'}s)^2$. We also have
$\kappa=2(n-1)\frac{r'}r$.

\medskip
{\bf Theorem 9.4. }{\it Let $\m$ be a complete QCH K\"ahler
manifold of dimension $2n\ge 6$ and an analytic Riemannian metric
$g$. If $\kappa\ne 0$  then $int\ B=\emptyset$, the set  $U=\{x\in
M:\kappa(x)\ne 0\}$ is open and dense in $M$ and every point $x\in
U$ has a neighborhood $V\subset U$ biholomorphic to the manifold
$(\a,\be)\times P$ with the metric (6.1) where $s\ne
0,f=\frac{2rr'}s, s\in\Bbb R$ with $r\in C^{\om}(\Bbb R)$, $P$ is
a circle bundle over a K\"ahler manifold $(N,h)$ and $h$ is a
K\"ahler metric of constant holomorphic sectional curvature  or
$b=0$ and  $M$ has a metric of constant holomorphic sectional
curvature. }
\medskip
{Proof.} The proof follows from Section 3, Lemma 9.1, Theorem 5.1
and from [D-M-2] (see Theorem 18.1).$\k$

\bigskip

{\it Acknowledgments.} The author would like to thank the referee
for his valuable remarks which improved the paper.
\bigskip
\centerline{\bf References.}
\par
\medskip
\cite{Bes}  A. L. Besse {\it Einstein manifolds}, Ergebnisse,
ser.3, vol. 10, Springer-Verlag, Berlin-Heidelberg-New York, 1987.
\par
\medskip
\cite{Ber}  L. B\'erard Bergery,{\it Sur de nouvelles vari\'et\'es
riemanniennes d'Einstein}, Pu\-bl. de l'Institute E. Cartan
(Nancy) {\bf 4},(1982), 1-60.
\par
\medskip
\cite{D-M-1} A. Derdzi\'nski, G. Maschler {\it Special
K\"ahler-Ricci potentials on compact K\" ahler manifolds}, J.
reine angew. Math. 593 (2006), 73-116.
 \par
\medskip
\cite{D-M-2}  A. Derdzi\'nski, G.  Maschler {\it Local
classification of conformally-Einstein  K\"ahler metrics in higher
dimensions}, Proc. London Math. Soc. (3) 87 (2003), no. 3,
779-819.

\par
\medskip
\cite{G-M-1} G.Ganchev, V. Mihova {\it K\"ahler manifolds of
quasi-constant holomorphic sectional curvatures}, Cent. Eur. J.
Math. 6(1),(2008), 43-75.
\par
\medskip
\cite{G-M-2} G.Ganchev, V. Mihova {\it Warped product K\"ahler
manifolds and Bochner-K\"ahler metrics}, J. Geom. Phys. 58(2008),
803-824.
\par
\medskip
\cite{H-H} F. Hirzebruch, H. Hopf  {\it Felder von
Fl\"achenelementen in 4-dimensionalen Mannigfaltigkeiten.}, Math.
Ann. 136, (1958), 156-172.

 \par
\medskip
\cite{K} S. Kobayashi {\it Principal fibre bundles with the
1-dimensional toroidal group} T\^ohoku Math.J. {\bf 8},(1956)
29-45.
\par
\medskip
\cite{K-N} S. Kobayashi and K. Nomizu {\it Foundations of
Differential Geometry}, vol.2, Interscience, New York  1963
\par
\medskip
\cite{ON} B. O'Neill, {\it The fundamental equations of a
submersion}, Mich. Math. J. {\bf 13},(1966), 459-469.
\par
\medskip
\cite{S} P. Sentenac ,{\it Construction d'une m\'etrique
d'Einstein  sur la somme de deux projectifs complexes de dimension
2}, G\'eom\'etrie riemannienne en dimension 4 ( S\'eminaire Arthur
Besse 1978-1979) Cedic-Fernand Nathan, Paris (1981), pp. 292-307.

\par
\medskip

\par
\medskip
Institute of Mathematics

Cracow University of Technology

Warszawska 24

31-155      Krak\'ow,  POLAND.

 E-mail address: wjelon\@pk.edu.pl
\bigskip

\enddocument